\documentclass[11pt,reqno]{amsart}
\usepackage{amsmath,amsthm,amssymb,verbatim,enumerate,ifthen}
\usepackage[mathscr]{eucal}
\usepackage[utf8]{inputenc}
\def\N{\mathbb{N}}
\def\R{\mathbb{R}}

\def\C{\mathscr{C}}

\def\CM{\mathscr{CM}}
\def\CP{\mathscr{CP}}

\def\muhat{\widehat{\mu}}

\newtheorem{theorem}{Theorem}
\newtheorem*{theorem*}{Theorem}
\def\Thm#1#2{\ifthenelse{\equal{#1}{*}}{\begin{theorem*}#2\end{theorem*}}
  {\begin{theorem}\label{T#1}#2\end{theorem}}}
\newtheorem{Atheorem}{Theorem}

\def\thm#1{Theorem~\ref{T#1}}

\newtheorem{proposition}[theorem]{Proposition}
\newtheorem*{proposition*}{Proposition}
\def\Prp#1#2{\ifthenelse{\equal{#1}{*}}{\begin{proposition*}#2\end{proposition*}}
             {\begin{proposition}\label{P#1}#2\end{proposition}}}
\def\prp#1{Proposition~\ref{P#1}}

\newtheorem{corollary}[theorem]{Corollary}
\newtheorem*{corollary*}{Corollary}
\def\Cor#1#2{\ifthenelse{\equal{#1}{*}}{\begin{corollary*}#2\end{corollary*}}
             {\begin{corollary}\label{C#1}#2\end{corollary}}}
\def\cor#1{Corollary~\ref{C#1}}

\newtheorem{lemma}[theorem]{Lemma}
\newtheorem*{lemma*}{Lemma}
\def\Lem#1#2{\ifthenelse{\equal{#1}{*}}{\begin{lemma*}#2\end{lemma*}}
             {\begin{lemma}\label{L#1}#2\end{lemma}}}
\def\lem#1{Lemma~\ref{L#1}}

\newtheorem{example}[theorem]{Example}
\newtheorem*{example*}{Example}
\def\Exa#1#2{\ifthenelse{\equal{#1}{*}}{\begin{example*}\rm #2\end{example*}}
             {\begin{example}\label{Ex#1}\rm #2\end{example}}}

\newtheorem{problem}[theorem]{Problem}

\theoremstyle{definition}
\newtheorem{definition}[theorem]{Definition}

\newtheorem{remark}[theorem]{Remark}
\newtheorem*{remark*}{Remark}
\def\Rem#1#2{\ifthenelse{\equal{#1}{*}}{\begin{remark*}\rm #2\end{remark*}}
             {\begin{remark}\label{R#1}\rm #2\end{remark}}}

\def\eq#1{{\rm(\ref{E#1})}}
\def\Eq#1#2{\ifthenelse{\equal{#1}{*}}
  {\begin{equation*}\begin{aligned}#2\end{aligned}\end{equation*}}
  {\begin{equation}\begin{aligned}\label{E#1}#2\end{aligned}\end{equation}}}

\begin{document}
\begin{flushright}
\end{flushright}
\vspace{5mm}

\date{\today}
\title{On the equality of two-variable general functional means}

\author[L. Losonczi]{L\'aszl\'o Losonczi}
\address[L. Losonczi]{Faculty of Economics, University of Debrecen, H-4028 Debrecen, B\"osz\"orm\'enyi \'ut 26, Hungary}
\email{laszlo.losonczi@econ.unideb.hu}

\author[Zs. P\'ales]{Zsolt P\'ales}
\address[Zs. P\'ales]{Institute of Mathematics, University of Debrecen, H-4002 Debrecen, Pf.\ 400, Hungary}
\email{pales@science.unideb.hu}

\author[A. Zakaria]{Amr Zakaria}
\address[A. Zakaria]{Doctoral School of Mathematical and Computational Sciences, University of Debrecen, H-4002 Debrecen, Pf.\ 400, Hungary; Department of Mathematics, Faculty of Education, Ain Shams University, Cairo 11341, Egypt}
\email{amr.zakaria@edu.asu.edu.eg}

\subjclass[2010]{39B22, 39B52}
\keywords{generalized quasiarithmetic mean; equality problem; functional equation}


\thanks{The research of the second author was supported by the EFOP-3.6.1-16-2016-00022 and the EFOP-3.6.2-16-2017-00015 projects. These projects are co-financed by the European Union and the European Social Fund.}

\begin{abstract}
Given two functions $f,g:I\to\R$ and a probability measure $\mu$ on the Borel subsets of $[0,1]$, the two-variable mean $M_{f,g;\mu}:I^2\to I$ is defined by
\Eq{*}{
   M_{f,g;\mu}(x,y)
      :=\bigg(\frac{f}{g}\bigg)^{-1}\left(
            \frac{\int_0^1 f\big(tx+(1-t)y\big)d\mu(t)}
                 {\int_0^1 g\big(tx+(1-t)y\big)d\mu(t)}\right)
   \qquad(x,y\in I).
}
This class of means includes quasiarithmetic as well as Cauchy and Bajraktarevi\'c means. The aim of this paper is, for a fixed probability measure $\mu$, to study their equality problem, i.e., to characterize those pairs of functions $(f,g)$ and $(F,G)$ such that
\Eq{*}{
   M_{f,g;\mu}(x,y)=M_{F,G;\mu}(x,y)   \qquad(x,y\in I)
}
holds. Under at most sixth-order differentiability assumptions for the unknown functions $f,g$ and $F,G$, we obtain several necessary conditions for the solutions of the above functional equation. For two particular measures, a complete description is obtained. These latter results offer eight equivalent conditions for the equality of Bajraktarevi\'c means and of Cauchy means. 
\end{abstract}

\maketitle

\section{\bf Introduction}

Throughout this paper $I$ will stand for a nonempty open real interval. In the sequel, the classes of continuous strictly monotone and continuous positive real-valued functions defined on $I$ will be denoted by $\CM(I)$ and $\CP(I)$, respectively.

In general, a continuous function $M:I^2\to I$ is called a
\textit{two-variable mean} on $I$ if the so-called mean value inequality
\Eq{MI}{
  \min(x,y)\leq M(x,y)\leq \max(x,y)   \qquad(x,y\in I)
}
holds. If, for $x\neq y$, both of the inequalities in \eq{MI} are strict, then $M$ is called a \textit{two-variable strict mean}. The arithmetic and geometric means are well known instances for strict means on $\R_+$.

Given a function $\varphi\in\CM(I)$, the \textit{two-variable quasiarithmetic mean} generated by $\varphi$ is the function $A_\varphi:I^2\to I$ defined by
\Eq{*}{
  A_\varphi(x,y)
      :=\varphi^{-1}\bigg(\frac{\varphi(x)+\varphi(y)}{2}\bigg)
   \qquad(x,y\in I).
}
The systematic treatment of these means was first given by Hardy, Littlewood and P\'olya \cite{HarLitPol34}. The most basic problem, the characterization of the equality of these means, states that
$A_\varphi$ and $A_\psi$ are equal to each other if and only if there exist two real constants $a\neq0$ and $b$ such that $\psi=a\varphi+b$.

The characterization of quasiarithmetic means was solved independently by Kolmogorov \cite{Kol30}, Nagumo \cite{Nag30}, de Finetti \cite{Def31} for the case when the number of variables is non-fixed. For the two-variable case, Acz\'el \cite{Acz47b}, \cite{Acz48a}, \cite{Acz48b}, \cite{Acz56c} (seel also \cite{AczDho89}), proved a characterization theorem involving the notion of bisymmetry. This result was extended to the $n$-variable case by Maksa--M\"unnich--Mokken \cite{MunMakMok99}. 

In this paper, we consider the following generalization of
quasiarithmetic means, which was introduced in \cite{LosPal08} and also investigated in \cite{LosPal11a}. Given two continuous functions $f,g:I\to\R$ with $g\in\CP(I)$, $f/g\in\CM(I)$ and a probability measure $\mu$ on the Borel subsets
of $[0,1]$, the two-variable mean $M_{f,g;\mu}:I^2\to I$ is defined by
\Eq{MM}{
   M_{f,g;\mu}(x,y)
      :=\bigg(\frac{f}{g}\bigg)^{-1}\left(
            \frac{\int_0^1 f\big(tx+(1-t)y\big)d\mu(t)}
                 {\int_0^1 g\big(tx+(1-t)y\big)d\mu(t)}\right)
   \qquad(x,y\in I).
}
Means of the above form, will be called \textit{generalized quasiarithmetic means}.

The first important particular case of this definition is when $\mu=\frac{\delta_0+\delta_1}{2}$. Here and in the sequel, $\delta_\tau$ will denote the Dirac measure concentrated at the point $\tau\in[0,1]$.
If $\varphi\in\CM(I)$, and $p\in\CP(I)$, then $M_{\varphi\cdot p,p;\mu}=B_{\varphi,p}$, which is defined by
\Eq{*}{
  B_{\varphi,p}(x,y)
      :=\varphi^{-1}\bigg(\frac{p(x)\varphi(x)+p(y)\varphi(y)}{p(x)+p(y)}\bigg)
   \qquad(x,y\in I)
}
and is called the \emph{two-variable Bajraktarevi\'c mean} (cf. Bajraktarevi\'c \cite{Baj58}, \cite{Baj63}). By taking $p=1$, we can see that this class of means includes quasiarithmetic means. Assuming 6 times continuous differentiability, the equality problem of these means was solved by Losonczi \cite{Los99}, \cite{Los06b}. 

The second important particular case is when $\mu$ is the Lebesgue measure on $[0,1]$ and $\varphi,\psi:I\to I$ are continuously differentiable functions with $\psi'\in\CP(I)$ and
$\varphi'/\psi'\in\CM(I)$. Then, by the fundamental theorem of Calculus, one can easily see that $M_{\varphi',\psi';\mu}=C_{\varphi,\psi}$, which is defined by
\Eq{*}{
  C_{\varphi,\psi}(x,y)
      :=\begin{cases}
        \bigg(\dfrac{\varphi'}{\psi'}\bigg)^{-1}
             \bigg(\dfrac{\varphi(y)-\varphi(x)}{\psi(y)-\psi(x)}\bigg)
            &\mbox{if }x\neq y\\
         x  &\mbox{if }x=y
        \end{cases}
   \qquad(x,y\in I).
}
Assuming 7 times continuous differentiability, the equality problem of these means was solved by Losonczi \cite{Los00a}.

The equality problem of means in various classes of two-variable means has been solved. We refer here to Losonczi's works \cite{Los99}, \cite{Los00a}, \cite{Los02a}, \cite{Los03a}, \cite{Los06b} where the equality of two-variable means is characterized. A key idea in these papers, under high order differentiablity assumptions, is to calculate and then compare the partial derivatives of the means at points of the form
$(x,x)$. A similar problem, the mixed equality problem of
quasiarithmetic and Lagrangian means was solved by P\'ales \cite{Pal11}.

The aim of this paper is to study the \textit{equality problem} of generalized quasiarithmetic means, i.e., to characterize those pairs of functions $(f,g)$ and $(F,G)$ such that
\Eq{P}{
   M_{f,g;\mu}(x,y)=M_{F,G;\mu}(x,y)   \qquad(x,y\in I)
}
holds. Due to the complexity of the problem, we will not solve it in its natural generality. In our final main results we consider the cases when the measure $\mu$ is either of the form $\frac{\delta_0+\delta_1}{2}$ or is the Lebesgue measure on $[0,1]$. For these two cases, we shall need sixth-order differentiability properties of the unknown functions $f,g$ and $F,G$.

\section{\bf Basic results}

Given a Borel probability measure $\mu$ on the interval $[0,1]$, we define the \emph{first moment} and the \emph{$n$th centralized moment} of $\mu$ by
\Eq{*}{
  \muhat_1:=\int_0^1 t d\mu(t) \qquad\mbox{and}\qquad
  \mu_n:=\int_0^1 (t-\muhat_1)^n d\mu(t) \qquad(n\in\N\cup\{0\}).
}
Clearly, $\mu_0=1$ and $\mu_1=0$. It is also obvious that
$\mu_{2n}\geq0$ and equality can hold if and only if $\mu$ is the Dirac measure $\delta_{\muhat_1}$.

In order to describe the regularity conditions related to the two unknown functions $f,g$ generating the mean $M_{f,g;\mu}$, we introduce some notations. The class $\C_0(I)$ consists of all those pairs $(f,g)$ of
continuous functions $f,g:I\to\R$ such that $g\in\CP(I)$ and $f/g\in\CM(I)$.
For $n\in\N$, we say that the pair $(f,g)$ is in the class $\C_n(I)$ if $f,g$ are
$n$-times continuously differentiable functions such that $g\in\CP(I)$ and the
function $f'g-fg'$ does not vanish anywhere on $I$. Obviously, this latter
condition implies that $f/g$ is strictly monotone, i.e., $f/g\in\CM(I)$.

For $(f,g)\in\C_2(I)$, we also introduce the notation
\Eq{*}{
  \Phi_{f,g}:=\frac{W_{f,g}^{2,0}}{W_{f,g}^{1,0}}\qquad\mbox{and}\qquad
  \Psi_{f,g}:=-\frac{W_{f,g}^{2,1}}{W_{f,g}^{1,0}},
}
where the $(i,j)$-order Wronskian operator $W^{i,j}$ is defined in terms of the $i$th and $j$th derivatives by
\Eq{*}{
  W_{f,g}^{i,j}:=\begin{vmatrix}f^{(i)}&f^{(j)}\\g^{(i)}&g^{(j)}\end{vmatrix}.
}

Our first result establishes a formula for the higher-order
derivatives of $f$ and $g$ as well as for their $(i,j)$-order Wronskian
in terms of the functions $\Phi_{f,g}$ and $\Psi_{f,g}$.

\Lem{0}{Let $(f,g)\in\C_n(I)$, where $n\geq2$ and define two sequences
$(\varphi_i)$ and $(\psi_i)$ by the recursions
\Eq{Rec}{
      \varphi_0&:=0,&\qquad
  \varphi_{i+1}&:=\varphi'_i+\varphi_i\Phi_{f,g}+\psi_i&
      \qquad (i&\in\{0,\dots,n-1\}),\\
         \psi_0&:=1,&
     \psi_{i+1}&:=\varphi_i\Psi_{f,g}+\psi'_i &
      \qquad (i&\in\{0,\dots,n-1\}).
}
Then, for $h\in\{f,g\}$,
\Eq{h}{
   h^{(i)}=\varphi_ih'+\psi_ih \qquad (i\in\{0,\dots,n\})
}
and
\Eq{W}{
  W_{f,g}^{i,j}=\begin{vmatrix}\varphi_i&\varphi_j\\\psi_i&\psi_j\end{vmatrix}
               \cdot W_{f,g}^{1,0}\qquad (i,j\in\{0,\dots,n\}).
}
In particular, 
\Eq{W01}{
  W_{f,g}^{i,0}=\varphi_i W_{f,g}^{1,0},\qquad
  W_{f,g}^{i,1}=-\psi_i W_{f,g}^{1,0}
  \qquad (i\in\{0,\dots,n\}).
}}

\begin{proof}Consider the second-order linear differential equation
\Eq{de1}{
  \begin{vmatrix}f''&f'&f\\g''&g'&g\\h''&h'&h\end{vmatrix}
   =W_{f,g}^{1,0}h''-W_{f,g}^{2,0}h'+W_{f,g}^{2,1}h=0
}
for the unknown function $h:I\to\R$. Obviously, \eq{de1} is satisfied for
$h\in\{f,g\}$. We can see that \eq{de1} is equivalent to the equation
\Eq{de}{
  h''=\Phi_{f,g}h'+\Psi_{f,g}h,
}
i.e., $h\in\{f,g\}$ is also a solution of \eq{de}. Observe that \eq{de} is the $i=2$ particular case of \eq{h}.)

The equality \eq{h} trivially holds if $i=0$. Assume that \eq{h} has been
proved for some $i\in\{0,\dots,n-1\}$. Then, using \eq{h}, \eq{de}, and
\eq{Rec}, we get
\Eq{*}{
 h^{(i+1)}=\big(h^{(i)}\big)'&=\big(\varphi_ih'+\psi_ih\big)'
   =\varphi_ih''+\varphi'_ih'+\psi_ih'+\psi'_ih\\
   &=\varphi_i(\Phi_{f,g}h'+\Psi_{f,g}h)+\varphi'_ih'+\psi_ih'+\psi'_ih
    =\varphi_{i+1}h'+\psi_{i+1}h,
}
which completes the proof of the induction.

The equality in \eq{W} follows from \eq{h}. Indeed, applying \eq{h} for $h=f$ and $h=g$, we obtain
\Eq{*}{
   W_{f,g}^{i,j}=\begin{vmatrix}f^{(i)}&f^{(j)}\\g^{(i)}&g^{(j)}\end{vmatrix}
   =\begin{vmatrix}\varphi_if'+\psi_if&\varphi_jf'+\psi_jf\\
                   \varphi_ig'+\psi_ig&\varphi_jg'+\psi_jg\end{vmatrix}
   =\begin{vmatrix}f'&f\\g'&g\end{vmatrix}\cdot
    \begin{vmatrix}\varphi_i&\varphi_j\\\psi_i&\psi_j\end{vmatrix}
   =\begin{vmatrix}\varphi_i&\varphi_j\\\psi_i&\psi_j\end{vmatrix}\cdot W_{f,g}^{1,0}.
}
Hence, the proof of the lemma is complete.
\end{proof}

In the sequel we shall need the first few members of the sequences $(\varphi_i)$ and $(\psi_i)$ constructed in \eq{Rec}. For the sake of convenience, we list, for small $i$, the first few members of them:
\Eq{*}{
   \varphi_1&=1, & \varphi_2&=\Phi_{f,g}, 
      & \varphi_3&=\Phi'_{f,g}+\Phi^2_{f,g}+\Psi_{f,g},
      & \varphi_4&=\Phi''_{f,g}+3\Phi'_{f,g}\Phi_{f,g}+\Phi^3_{f,g}
          +2\Phi_{f,g}\Psi_{f,g}+2\Psi'_{f,g},\\
   \psi_1&=0, &\psi_2&=\Psi_{f,g},
      & \psi_3&=\Phi_{f,g}\Psi_{f,g}+\Psi'_{f,g},
      & \psi_4&=\big(2\Phi'_{f,g}+\Phi^2_{f,g}\big)\Psi_{f,g}
          +\Phi_{f,g}\Psi'_{f,g}+\Psi^2_{f,g}+\Psi''_{f,g},
}
\Eq{*}{
   \varphi_5&=\Phi'''_{f,g}+4\Phi''_{f,g}\Phi_{f,g}+3{\Phi'}^2_{f,g}
       +6\Phi'_{f,g}\Phi^2_{f,g}+\Phi^4_{f,g}+\big(4\Phi'_{f,g}+3\Phi^2_{f,g}\big)\Psi_{f,g}
       +5\Phi_{f,g}\Psi'_{f,g}+\Psi^2_{f,g}+3\Psi''_{f,g},& \\
   \psi_5&=2\Phi_{f,g}\Psi^2_{f,g}
          +\big(3\Phi''_{f,g}+5\Phi'_{f,g}\Phi_{f,g}+\Phi^3_{f,g}\big)\Psi_{f,g}
          +\big(3\Phi'_{f,g}+\Phi^2_{f,g}\big)\Psi'_{f,g}+\Phi_{f,g}\Psi''_{f,g}
          +4\Psi'_{f,g}\Psi_{f,g}+\Psi'''_{f,g},
}
\Eq{*}{
   \varphi_6&=\Phi''''_{f,g}+5\Phi'''_{f,g}\Phi_{f,g}+10\Phi''_{f,g}\Phi'_{f,g}
        +10\Phi''_{f,g}\Phi^2_{f,g}+10\Phi'_{f,g}\Phi^3_{f,g}
        +15{\Phi'}^2_{f,g}\Phi_{f,g}+\Phi^5_{f,g}+3\Phi_{f,g}\Psi^2_{f,g}\\
        &\quad+\big(7\Phi''_{f,g}+15\Phi'_{f,g}\Phi_{f,g}+4\Phi^3_{f,g}\big)\Psi_{f,g}
         +\big(12\Phi'_{f,g}+9\Phi^2_{f,g}\big)\Psi'_{f,g}+9\Phi_{f,g}\Psi''_{f,g}
         +6\Psi'_{f,g}\Psi_{f,g}+4\Psi'''_{f,g},\\
     \psi_6&=\big(6\Phi'_{f,g}+3\Phi^2_{f,g}\big)\Psi^2_{f,g}
        +\big(4\Phi'''_{f,g}+9\Phi''_{f,g}\Phi_{f,g}+8{\Phi'}^2_{f,g}
        +9\Phi'_{f,g}\Phi^2_{f,g}+\Phi^4_{f,g}\big)\Psi_{f,g}
        +\big(4\Phi'_{f,g}+\Phi^2_{f,g}\big)\Psi''_{f,g}\\
        &\quad+\big(6\Phi''_{f,g}+7\Phi'_{f,g}\Phi_{f,g}+\Phi^3_{f,g}\big)\Psi'_{f,g}
        +\Phi_{f,g}\big(\Psi'''_{f,g}+9\Psi'_{f,g}\Psi_{f,g}\big)
        +\Psi^3_{f,g}+4{\Psi'}^2_{f,g}+7\Psi''_{f,g}\Psi_{f,g}+\Psi''''_{f,g}.
}

We say that \textit{two pairs of functions $(f,g),(F,G)\in\C_0(I)$ are equivalent}, denoted by $(f,g)\sim(F,G)$, if there exists a nonsingular
$2\times2$-matrix $A$ (with real entries) such that
\Eq{A}{
   \binom{F}{G}=A\binom{f}{g}.
}
In other words, $(f,g)\sim(F,G)$ holds if there exist four real constants $a,b,c,d$ with $ad\neq bc$ such that $F=af+bg$ and $G=cf+dg$.

The remaining auxiliary results of this section were obtained in \cite{LosPal08} and \cite{PalZak20a}. The property of equivalence in the class $\C_2(I)$ is completely characterized by the following result.

\Lem{3}{{\rm(\cite[Theorem 2.1]{PalZak20a})} Let $(f,g),(F,G)\in\C_2(I)$. Then $(f,g)\sim(F,G)$ holds if and only if
\Eq{3}{
  \Phi_{f,g}=\Phi_{F,G}\qquad\mbox{and}\qquad \Psi_{f,g}=\Psi_{F,G}.
}}

The next result characterizes the mean $M_{f,g;\mu}$ via an implicit equation.

\Lem{1}{{\rm(\cite[Lemma 1]{LosPal08}, \cite[Lemma 1.1]{PalZak20a})} Let $(f,g)\in\C_0(I)$ and $\mu$ be a Borel probability measure on $[0,1]$.
Then for all $x,y\in I$, the value $z=M_{f,g;\mu}(x,y)$ is the unique solution
of the equation
\Eq{1}{
  \int_0^1 \begin{vmatrix}f\big(tx+(1-t)y\big) & f(z)\\
            g\big(tx+(1-t)y\big) & g(z)\end{vmatrix}d\mu(t)=0.
}}

As a consequence, the next lemma shows that the equivalent pairs of generating functions determine identical means.

\Lem{2}{{\rm(\cite{LosPal08}, \cite{PalZak20a})} Let $(f,g),(F,G)\in\C_0(I)$ and $\mu$ be a Borel probability measure on $[0,1]$. Assume that $(f,g)\sim(F,G)$. Then $M_{F,G;\mu}=M_{f,g;\mu}$.}

\section{\bf Higher-order directional derivatives of generalized quasiarithmetic means}

\Lem{n}{Let $n\in\N$, $(f,g)\in\C_n(I)$ and $\mu$ be a Borel probability measure on $[0,1]$. Then $M_{f,g;\mu}$ is $n$-times continuously differentiable on $I\times I$.}

\begin{proof} The proof of this statement requires the use of standard calculus rules and a standard argument. One can verify that the inverse of $f/g$ and the maps $(x,y)\mapsto\int_0^1 h(tx+(1-t)y)d\mu(t)$ (where $h\in\{f,g\}$) are $n$-times differentiable on $(f/g)(I)$ and $I^2$, respectively. 
\end{proof}

In what follows, we deduce explicit formulae for the high-order directional
derivatives of $M_{f,g;\mu}$ at the diagonal points of the Cartesian product
$I\times I$. Given $(f,g)\in\C_0(I)$ and a fixed element $x\in I$, define 
the function $m_x=m_{x;f,g;\mu}$ in a neigborhood of zero by
\Eq{m}{
  m_x(u):=m_{x;f,g;\mu}(u):=M_{f,g;\mu}\big(x+(1-\hat{\mu}_1)u,x-\hat{\mu}_1 u\big),
}
where $\hat{\mu}_1$ denotes the first moment of the measure $\mu$.

\Prp{der}{Let $n\in\N$, $(f,g)\in\C_n(I)$, and $\mu$ be a Borel 
probability measure on $[0,1]$. Then, for fixed $x\in I$, the function $m_x$ defined by \eq{m} is $n$-times continuously differentiable at the origin and
\Eq{der}{
  \sum_{i=0}^n\binom{n}{i}\mu_i\begin{vmatrix}f^{(i)}(x) & (f\circ m_x)^{(n-i)}(0)\\
    g^{(i)}(x) & (g\circ m_x)^{(n-i)}(0)\end{vmatrix}=0.
}
Furthermore, $m_x(0)=x$ and in the cases $n=1,2,3,4,5,6$, we have
\Eq{*}{
       m_x'(0)&=0,\\
      m_x''(0)&=\mu_2\Phi_{f,g}(x),\\
     m_x'''(0)&=\mu_3\big(\Phi'_{f,g}+\Phi^2_{f,g}+\Psi_{f,g}\big)(x),\\
    m_x''''(0)&=-3\mu_2^2\big(\Phi_{f,g}^3+2\Phi_{f,g}\Psi_{f,g}\big)(x)+\mu_4\big(\Phi''_{f,g}+3\Phi'_{f,g}\Phi_{f,g}+\Phi^3_{f,g} +2\Phi_{f,g}\Psi_{f,g}+2\Psi'_{f,g}\big)(x),\\
   m_x'''''(0)&=-10\mu_3\mu_2\big(\Phi_{f,g}^2\Phi'_{f,g}+\Phi_{f,g}^4
      +\big(\Phi'_{f,g}+3\Phi_{f,g}^2\big)\Psi_{f,g}
      +\Phi_{f,g}\Psi'_{f,g}+\Psi_{f,g}^2\big)(x)\\
       &\quad+\mu_5\big(\Phi'''_{f,g}+4\Phi''_{f,g}\Phi_{f,g}+3{\Phi'}^2_{f,g}
       +6\Phi'_{f,g}\Phi^2_{f,g}+\Phi^4_{f,g}+\\
       &\qquad\qquad\qquad+\big(4\Phi'_{f,g}+3\Phi^2_{f,g}\big)\Psi_{f,g}
       +5\Phi_{f,g}\Psi'_{f,g}+\Psi^2_{f,g}+3\Psi''_{f,g}\big)(x),
}
and
\Eq{*}{
  m_x''''''(0)&=+15\mu_2^3\big(-\Phi'_{f,g}\Phi^3_{f,g}+2\Phi_{f,g}^5
      +8\Phi^3_{f,g}\Psi_{f,g}+6\Phi_{f,g}\Psi^2_{f,g}\big)(x)\\
    &\quad-10\mu_3^2(
    {\Phi'}^2_{f,g}\Phi_{f,g}+2\Phi'_{f,g}\Phi_{f,g}^3+\Phi_{f,g}^5+3\Phi_{f,g}\Psi_{f,g}^2
    +4\big(\Phi'_{f,g}\Phi_{f,g}+\Phi_{f,g}^3\big)\Psi_{f,g}\\
    &\qquad\qquad\qquad+2\big(\Phi'_{f,g}+\Phi_{f,g}^2\big)\Psi'_{f,g}
      +2\Psi'_{f,g}\Psi_{f,g}\big)(x)\\
    &\quad-15\mu_2\mu_4\big(\Phi''_{f,g}\Phi^2_{f,g}+3\Phi'_{f,g}\Phi^3_{f,g}
      +\Phi^5_{f,g}+3\Phi_{f,g}\Psi_{f,g}^2
      +\big(\Phi''_{f,g}+5\Phi'_{f,g}\Phi_{f,g}+4\Phi^3_{f,g}\big)\Psi_{f,g}\\
    &\qquad\qquad\qquad+3\Phi^2_{f,g}\Psi'_{f,g}+\Phi_{f,g}\Psi''_{f,g}
      +2\Psi'_{f,g}\Psi_{f,g}\big)(x)\\
    &\quad+\mu_6\big(\Phi''''_{f,g}+5\Phi'''_{f,g}\Phi_{f,g}+10\Phi''_{f,g}\Phi'_{f,g}
      +10\Phi''_{f,g}\Phi^2_{f,g}+10\Phi'_{f,g}\Phi^3_{f,g}
      +15{\Phi'}^2_{f,g}\Phi_{f,g}+\Phi^5_{f,g}\\
    &\qquad\qquad\qquad+3\Phi_{f,g}\Psi^2_{f,g}+\big(7\Phi''_{f,g}+15\Phi'_{f,g}\Phi_{f,g}
      +4\Phi^3_{f,g}\big)\Psi_{f,g}+\big(12\Phi'_{f,g}+9\Phi^2_{f,g}\big)\Psi'_{f,g}\\
    &\qquad\qquad\qquad+9\Phi_{f,g}\Psi''_{f,g}+6\Psi'_{f,g}\Psi_{f,g}+4\Psi'''_{f,g}\big)(x).
}}

\begin{proof} Let $x\in I$ be fixed. Then the $n$-times continuously differentiability of $m_x$ at the origin is a consequence of \lem{n}. By \lem{1}, for small $u$, we have that $u$ is the solution of the equation
\Eq{*}{
  \int_0^1 \begin{vmatrix}f\big(x+(t-\hat{\mu}_1)u\big) & (f\circ m_x)(u)\\
            g\big(x+(t-\hat{\mu}_1)u\big) & (g\circ m_x)(u)\end{vmatrix}d\mu(t)=0.
}
Differentiating this equality $n$-times with respect to the variable $u$ and
using Leibniz's rule, we obtain
\Eq{*}{
  \sum_{i=0}^n\binom{n}{i}\int_0^1
  \begin{vmatrix}f^{(i)}\big(x+(t-\hat{\mu}_1)u\big) & (f\circ m_x)^{(n-i)}(u)\\
   g^{(i)}\big(x+(t-\hat{\mu}_1)u\big) & (g\circ m_x)^{(n-i)}(u)\end{vmatrix}
  (t-\hat{\mu}_1)^id\mu(t)=0.
}
Now the substitution $u=0$ implies \eq{der}.

The equality $m_x(0)=x$ follows from the definition of $m_x$. By $\mu_0=1$,
$\mu_1=0$, in the case $n=1$, \eq{der} reduces to
\Eq{*}{
  0=\begin{vmatrix}f(x) & (f\circ m_x)'(0)\\g(x) & (g\circ m_x)'(0)\end{vmatrix}
   =\begin{vmatrix}f(x) & f'(x)\\g(x) & g'(x)\end{vmatrix}m_x'(0),
}
which yields $m_x'(0)=0$.

To elaborate the condition \eq{der} in the cases $n=2,3,4,5,6$, we shall need
the following computational rules for the functions $h\in\{f,g\}$,
\Eq{*}{
  (h\circ m_x)'&=(h'\circ m_x)m_x',\\
  (h\circ m_x)''&=(h''\circ m_x)(m_x')^2+(h'\circ m_x)m_x'',\\
  (h\circ m_x)'''&=(h'''\circ m_x)(m_x')^3+3(h''\circ m_x)m_x''m_x'+(h'\circ m_x)m_x''',\\
  (h\circ m_x)''''&=(h''''\circ m_x)(m_x')^4+6(h'''\circ m_x)m_x''(m_x')^2
          +(h''\circ m_x)\big(4m_x'''m_x'+3(m_x'')^2\big)+(h'\circ m_x)m_x'''',\\
  (h\circ m_x)'''''&=(h'''''\circ m_x)(m_x')^5+10(h''''\circ m_x)m_x''(m_x')^3
          +(h'''\circ m_x)\big(10m_x'''(m_x')^2+15(m_x'')^2m_x'\big)\\
          &\qquad+(h''\circ m_x)\big(5m_x''''m_x'+10m_x'''m_x''\big)+(h'\circ m_x)m_x''''',\\
  (h\circ m_x)''''''&=(h''''''\circ m_x)(m_x')^6+15(h'''''\circ m_x)m_x''(m_x')^4
          +(h''''\circ m_x)\big(20m_x'''(m_x')^3+45(m_x'')^2(m_x')^2\big)\\
          &\qquad+(h'''\circ m_x)\big(15m_x''''(m_x')^2+60m_x'''m_x''m_x'+15(m_x'')^3\big)\\
          &\qquad+(h''\circ m_x)\big(6m_x'''''m_x'+15m_x''''m_x''+10(m_x''')^2\big)
            +(h'\circ m_x)m_x''''''.
}
Hence, using $m_x(0)=x$ and $m_x'(0)=0$, it follows that
\Eq{ids}{
  (h\circ m_x)'(0)&=0,\\
  (h\circ m_x)''(0)&=h'(x)m_x''(0),\\
  (h\circ m_x)'''(0)&=h'(x)m_x'''(0),\\
  (h\circ m_x)''''(0)&=h'(x)m_x''''(0)+3h''(x)(m_x''(0))^2,\\
  (h\circ m_x)'''''(0)&=h'(x)m_x'''''(0)+10h''(x)m_x'''(0)m_x''(0),\\
  (h\circ m_x)''''''(0)&=h'(x)m_x''''''(0)+h''(x)\big(15m_x''''(0)m_x''(0)+10(m_x'''(0))^2\big)
                      +15h'''(x)(m_x''(0))^3.
}

In the case $n=2$, using $\mu_0=1$ and $\mu_1=0$, equation \eq{der} yields
\Eq{*}{
   \begin{vmatrix}f(x) & (f\circ m)''(0)\\
    g(x) & (g\circ m_x)''(0)\end{vmatrix}
  +\mu_2\begin{vmatrix}f''(x) & (f\circ m)(0)\\
    g''(x) & (g\circ m_x)(0)\end{vmatrix} =0,
}
which, in view of the identities \eq{ids}, reduces to the equality
\Eq{*}{
  -W_{f,g}^{1,0}(x)m_x''(0)+\mu_2W_{f,g}^{2,0}(x)=0,
}
proving $m_x''(0)=\mu_2\Phi_{f,g}(x)$.

In the case $n=3$, using $\mu_0=1$, $\mu_1=0$, and $m'(0)=0$, equation
\eq{der} gives
\Eq{*}{
  \begin{vmatrix}f(x) & (f\circ m_x)'''(0)\\
    g(x) & (g\circ m_x)'''(0)\end{vmatrix}
  +\mu_3\begin{vmatrix}f'''(x) & (f\circ m_x)(0)\\
    g'''(x) & (g\circ m_x)(0)\end{vmatrix} =0,
}
which, applying the third identity in \eq{ids} for $h=f$ and $h=g$, implies
\Eq{*}{
  -W_{f,g}^{1,0}(x)m_x'''(0)+\mu_3W_{f,g}^{3,0}(x)=0.
}
Hence, using \eq{W01} for $i=3$, we get
\Eq{*}{
  m_x'''(0)=\mu_3\frac{W_{f,g}^{3,0}}{W_{f,g}^{1,0}}(x)
    =\mu_3\varphi_3(x)=\mu_3\big(\Phi'_{f,g}+\Phi^2_{f,g}+\Psi_{f,g}\big)(x).
}

In the case $n=4$, using $\mu_0=1$, $\mu_1=0$ and $m'(0)=0$, equation \eq{der} results
\Eq{*}{
  \begin{vmatrix}f(x) & (f\circ m_x)''''(0)\\
    g(x) & (g\circ m_x)''''(0)\end{vmatrix}
    +6\mu_2\begin{vmatrix}f''(x) & (f\circ m_x)''(0)\\
    g''(x) & (g\circ m_x)''(0)\end{vmatrix}
  +\mu_4\begin{vmatrix}f''''(x) & (f\circ m_x)(0)\\
    g''''(x) & (g\circ m_x)(0)\end{vmatrix} =0,
}
which, by the second and fourth identities in \eq{ids}, implies
\Eq{*}{
  -W_{f,g}^{1,0}(x)m_x''''(0)-3W_{f,g}^{2,0}(x)(m_x''(0))^2
    +6\mu_2W_{f,g}^{2,1}(x)m_x''(0)+\mu_4W_{f,g}^{4,0}(x)=0.
}
Thus, applying \eq{W01} and the formulae for $m_x''(0)$, $\varphi_2,\varphi_4$,
and $\psi_2$, we obtain
\Eq{*}{
  m_x''''(0)
    &=-3\varphi_2(x)(m_x''(0))^2-6\mu_2\psi_2(x)m_x''(0)+\mu_4\varphi_4(x)\\
    &=-3\mu_2^2\big(\Phi_{f,g}^3+2\Phi_{f,g}\Psi_{f,g}\big)(x)
      +\mu_4\big(\Phi''_{f,g}+3\Phi'_{f,g}\Phi_{f,g}+\Phi^3_{f,g}+2\Phi_{f,g}\Psi_{f,g}+2\Psi'_{f,g}\big)(x).
}

In the case $n=5$, using $\mu_0=1$, $\mu_1=0$ and $m'(0)=0$, equation \eq{der} results
\Eq{*}{
  \begin{vmatrix}f(x) & (f\circ m_x)'''''(0)\\
    g(x) & (g\circ m_x)'''''(0)\end{vmatrix}
    &+10\mu_2\begin{vmatrix}f''(x) & (f\circ m_x)'''(0)\\
    g''(x) & (g\circ m_x)'''(0)\end{vmatrix}\\
    &+10\mu_3\begin{vmatrix}f'''(x) & (f\circ m_x)''(0)\\
    g'''(x) & (g\circ m_x)''(0)\end{vmatrix}
  +\mu_5\begin{vmatrix}f'''''(x) & (f\circ m_x)(0)\\
    g'''''(x) & (g\circ m_x)(0)\end{vmatrix} =0,
}
which, by the identities \eq{ids}, implies
\Eq{*}{
  -W_{f,g}^{1,0}(x)m_x'''''(0)&-10W_{f,g}^{2,0}(x)m_x'''(0)m_x''(0)
    +10\mu_2W_{f,g}^{2,1}(x)m_x'''(0)\\&+10\mu_3W_{f,g}^{3,1}(x)m_x''(0)+\mu_5W_{f,g}^{5,0}(x)=0.
}
Thus, applying \eq{W01}, the above equality yields that
\Eq{*}{
  m_x'''''(0)=-10\varphi_2(x)m_x'''(0)m_x''(0)-10\mu_2\psi_2(x)m_x'''(0)
              -10\mu_3\psi_3(x)m_x''(0)+\mu_5\varphi_5(x).
}
Using the explicit formulae for $m_x''(0)$, $m_x'''(0)$,
$\varphi_2,\varphi_5$, and $\psi_2,\psi_3$,  
this expression simplifies to the statement.

Finally, in the case $n=6$, using $\mu_0=1$, $\mu_1=0$ and $m'(0)=0$, equation \eq{der} results
\Eq{*}{
  \begin{vmatrix}f(x) & (f\circ m_x)''''''(0)\\
    g(x) & (g\circ m_x)''''''(0)\end{vmatrix}
    &+15\mu_2\begin{vmatrix}f''(x) & (f\circ m_x)''''(0)\\
    g''(x) & (g\circ m_x)''''(0)\end{vmatrix}
     +20\mu_3\begin{vmatrix}f'''(x) & (f\circ m_x)'''(0)\\
    g'''(x) & (g\circ m_x)'''(0)\end{vmatrix}\\
    &+15\mu_4\begin{vmatrix}f''''(x) & (f\circ m_x)''(0)\\
    g''''(x) & (g\circ m_x)''(0)\end{vmatrix}
  +\mu_6\begin{vmatrix}f''''''(x) & (f\circ m_x)(0)\\
    g''''''(x) & (g\circ m_x)(0)\end{vmatrix} =0,
}
which, using the identities \eq{ids}, implies
\Eq{*}{
  &-W_{f,g}^{1,0}(x)m_x''''''(0)-W_{f,g}^{2,0}(x)\big(15m_x''''(0)m_x''(0)
    +10(m_x'''(0))^2\big)-15W_{f,g}^{3,0}(x)(m_x''(0))^3\\
    &+15\mu_2W_{f,g}^{2,1}(x)m_x''''(0)+20\mu_3W_{f,g}^{3,1}(x)m_x'''(0)
    +15\mu_4W_{f,g}^{4,1}(x)m_x''(0)+\mu_6W_{f,g}^{6,0}(x)=0.
}
Thus, applying \eq{W01}, we get
\Eq{*}{
  m_x''''''(0)
    &=-\varphi_2(x)\big(15m_x''''(0)m_x''(0)+10(m_x'''(0))^2\big)-15\varphi_3(x)(m_x''(0))^3\\
    &\quad-15\mu_2\psi_2(x)m_x''''(0)-20\mu_3\psi_3(x)m_x'''(0)
        -15\mu_4\psi_4(x)m_x''(0)+\mu_6\varphi_6(x).
}
Using the explicit formulae for $m_x''(0),m_x'''(0),m_x''''(0)$, 
$\varphi_2,\varphi_3,\varphi_6$, and $\psi_2,\psi_3,\psi_4$, 
this formula simplifies to the statement.
\end{proof}

\section{\bf Necessary and sufficient conditions for the equality of generalized quasiarithmetic means}

In what follows, given $(f,g),(F,G)\in\C_0(I)$ and a probability measure $\mu$ on $[0,1]$, we say that $M_{f,g;\mu}$ equals $M_{F,G;\mu}$ if they coincide at every point of $I^2$. We say that these two means are \emph{equal near the diagonal $\Delta(I):=\{(x,x)\mid x\in I\}$ of $I^2$}, if there exists an open set $U\subseteq I^2$ containing a dense subset $D$ of $\Delta(I)$ such that the two means are equal at every point of $U$.

\Lem{N0}{
Let $\mu$ be a Borel probability measure on $[0,1]$, let $n\in\N$ and let $(f,g),(F,G)\in\C_n(I)$. If $M_{f,g;\mu}$ equals $M_{F,G;\mu}$ near the diagonal of $I^2$, then, for all $k\in\{1,\dots,n\}$ and $x\in I$,
\Eq{mk}{
m_{x;f,g;\mu}^{(k)}(0)=m_{x;F,G;\mu}^{(k)}(0).
}
}

\begin{proof}
Let $U\subseteq I^2$ be an open set containing a dense subset $D$ of $\Delta(I)$ such that the two means are equal at every point of $U$. Let $x\in I$ be fixed such that $(x,x)\in D$. Define 
\Eq{*}{
  U_x:=\{u\in\R\mid (x+(1-\hat{\mu}_1)u,x-\hat{\mu}_1 u)\in U\}.
}
Then $U_x$ is a neighbourhood of $0$ (because $U$ is open), and the equality of the means on $U$ implies that, for any $u\in U_x$, 
\Eq{*}{
m_{x;f,g;\mu}(u)=m_{x;F,G;\mu}(u).
}
Therefore \eq{mk} holds for all $k\in\{1,\dots,n\}$ and $x\in I$ with $(x,x)\in D$.
Using the continuity, the density of $D$ yields that this equality holds for all $x\in I$. 
\end{proof}

In the subsequent lemmas, we will analyze the consequences of the equalities \eq{mk} for $k\in\{2,3,4,5,6\}$.

\Lem{N1}{
Let $\mu$ be a Borel probability measure on $[0,1]$ with $\mu_2>0$ and let $(f,g),(F,G)\in\C_2(I)$. If, for all $x\in I$,
\Eq{2d}{
m_{x;f,g;\mu}''(0)=m_{x;F,G;\mu}''(0)
}
holds, then
\Eq{Phi}{
\Phi_{f,g}=\Phi_{F,G}.
}
}

\begin{proof}
Applying the second-order formula of \prp{der}, the equality \eq{2d} implies 
\Eq{*}{
\mu_2\Phi_{f,g}=\mu_2\Phi_{F,G},
}
which, using $\mu_2>0$, proves \eq{Phi}.
\end{proof}

\Lem{N1.5}{
Let $\mu$ be a Borel probability measure on $[0,1]$ with $\mu_3\neq0$ and let $(f,g),(F,G)\in\C_3(I)$. If, for all $x\in I$,
\Eq{2+3d}{
m_{x;f,g;\mu}''(0)=m_{x;F,G;\mu}''(0) \qquad\mbox{and}\qquad
m_{x;f,g;\mu}'''(0)=m_{x;F,G;\mu}'''(0)
}
hold, then
\Eq{Phi+Psi}{
\Phi_{f,g}=\Phi_{F,G} \qquad\mbox{and}\qquad \Psi_{f,g}=\Psi_{F,G}.
}
}

\begin{proof}
The condition $\mu_3\neq0$ implies that $\mu_2>0$ is also valid.
Thus, the first equality in \eq{Phi+Psi} is a consequence of \lem{N1}.
Applying the third-order formula of \prp{der}, the second equality \eq{2+3d} implies 
\Eq{*}{
\mu_3\big(\Phi'_{f,g}+\Phi^2_{f,g}+\Psi_{f,g}\big)
=\mu_3\big(\Phi'_{F,G}+\Phi^2_{F,G}+\Psi_{F,G}\big),
}
which, using $\mu_3\neq0$ and the equality $\Phi_{f,g}=\Phi_{F,G}$, proves the last equation in \eq{Phi+Psi}.
\end{proof}

\Thm{N1.5}{
Let $\mu$ be a Borel probability measure on $[0,1]$ with $\mu_3\neq0$ and let $(f,g),(F,G)\in\C_3(I)$. Then the following assertions are equivalent:
\begin{enumerate}[(i)]
 \item The means $M_{f,g;\mu}$ and $M_{F,G;\mu}$ are equal on $I^2$.
 \item The means $M_{f,g;\mu}$ and $M_{F,G;\mu}$ are equal near the diagonal of  $I^2$.
 \item For all $x\in I$, the equalities in \eq{2+3d} hold.
 \item The equalities $\Phi_{f,g}=\Phi_{F,G}$ and $\Psi_{f,g}=\Psi_{F,G}$ hold on $I$.
 \item $(f,g)\sim(F,G)$ holds.
\end{enumerate}
}

\begin{proof}
The implication (i)$\Rightarrow$(ii) is obvious. The implication (ii)$\Rightarrow$(iii) is a consequence of \lem{N0}. The implication (iii)$\Rightarrow$(iv) follows from \lem{N1.5}. Finally, the implications (iv)$\Rightarrow$(v) and (v)$\Rightarrow$(i) are consequences of \lem{3} and \lem{2}, respectively.
\end{proof}

\Lem{N2}{
Let $\mu$ be a Borel probability measure on $[0,1]$ with $\mu_2>0$ and let $(f,g),(F,G)\in\C_4(I)$. If, for all $x\in I$,
\Eq{2+4d}{
m_{x;f,g;\mu}''(0)=m_{x;F,G;\mu}''(0) \qquad\mbox{and}\qquad
m_{x;f,g;\mu}''''(0)=m_{x;F,G;\mu}''''(0)
}
hold, then the equality $\Phi_{f,g}=\Phi_{F,G}=:\Phi$ holds on $I$ and there exists a constant $\gamma\in\R$ such that 
\Eq{R}{
\Psi_{f,g}-\Psi_{F,G}=2\gamma\big|W_{f,g}^{1,0}\big|^p, 
}
where 
\Eq{p}{
  p:=3\dfrac{\mu_2^2}{\mu_4}-1.
}}

\begin{proof}
The condition $\mu_2>0$ implies that $\mu_4>0$ is also valid. It follows from \lem{N1} and the first condition in \eq{2+4d} that $\Phi_{f,g}=\Phi_{F,G}$ holds on $I$. Using the formula for the fourth-order derivative by \prp{der}, the second equality in \eq{2+4d} simplifies to
\Eq{*}{
-6\mu_2^2\Phi\Psi_{f,g}+2\mu_4(\Phi\Psi_{f,g}+\Psi_{f,g}')
=-6\mu_2^2\Phi\Psi_{F,G}+2\mu_4(\Phi\Psi_{F,G}+\Psi_{F,G}'),
}
Therefore, we get the following first-order homogeneous linear differential equation for the difference function $R:=\Psi_{f,g}-\Psi_{F,G}$:
\Eq{R'}{
R'=\Big(3\frac{\mu_2^2}{\mu_4}-1\Big)\Phi R
=p\Phi R
=p\frac{W_{f,g}^{2,0}}{W_{f,g}^{1,0}} R
=p\frac{\big(W_{f,g}^{1,0}\big)'}{W_{f,g}^{1,0}} R.
}
The solution of this differential equation implies \eq{R} for some $\gamma\in\R$.
\end{proof}

\Lem{N2.5}{
Let $\mu$ be a Borel probability measure on $[0,1]$ with $\mu_3=0\neq\mu_5$ and let $(f,g),(F,G)\in\C_5(I)$. If, for all $x\in I$,
\Eq{2+4+5d}{
m_{x;f,g;\mu}''(0)=m_{x;F,G;\mu}''(0), \qquad
m_{x;f,g;\mu}''''(0)=m_{x;F,G;\mu}''''(0) \qquad\mbox{and}\qquad
m_{x;f,g;\mu}'''''(0)=m_{x;F,G;\mu}'''''(0)
}
hold, then the equality $\Phi_{f,g}=\Phi_{F,G}=:\Phi$ holds on $I$ and either $\Psi_{f,g}=\Psi_{F,G}$ holds on $I$ or there exists a nonzero real constant $\gamma$ such that
\Eq{S}{
 \Psi_{f,g}
  &=\gamma\big|W_{f,g}^{1,0}\big|^p-\frac{4+3p}2\Phi'-\frac{3+5p+3p^2}2\Phi^2,\\
  \Psi_{F,G}
  &=-\gamma\big|W_{f,g}^{1,0}\big|^p-\frac{4+3p}2\Phi'-\frac{3+5p+3p^2}2\Phi^2,
}
where $p$ is defined by \eq{p}.}

\begin{proof}
The condition $\mu_5\neq0$ implies that $\mu_2\mu_4\neq0$ is also valid. Therefore, by \lem{N1} and \lem{N2}, the first two conditions in \eq{2+4+5d} yield that $\Phi_{f,g}=\Phi_{F,G}$ and, with the notation $R:=\Psi_{f,g}-\Psi_{F,G}$, the equalities in \eq{R} and \eq{R'} hold for some $\gamma\in\R$. Then
\Eq{R''}{
  R''=(p\Phi R)'=p\Phi'R+p\Phi R'=(p\Phi'+p^2\Phi^2)R.
}
From \lem{N2} it follows that $R$ is either identically zero or nowhere zero in $I$. In the first case, we have that $\Psi_{f,g}=\Psi_{F,G}$. Thus, in the rest of the proof, we may assume that $R$ is nowhere zero in $I$, i.e., $\gamma\neq0$.

Using the formula for the fifth-order derivative by \prp{der}, the third condition of \eq{2+4+5d} simplifies to
\Eq{*}{
  (4\Phi'+3\Phi^2\big)\Psi_{f,g}
       +5\Phi\Psi'_{f,g}+\Psi^2_{f,g}+3\Psi''_{f,g}
  =  (4\Phi'+3\Phi^2\big)\Psi_{F,G}
       +5\Phi\Psi'_{F,G}+\Psi^2_{F,G}+3\Psi''_{F,G}
}
Define $S:=\Psi_{f,g}+\Psi_{F,G}$. Then the above equality yields
\Eq{*}{
  (4\Phi'+3\Phi^2\big)R+5\Phi R'+RS+3R''=0.
}
Using \eq{R''}, we finally get that
\Eq{*}{
  \big((4+3p)\Phi'+(3+5p+3p^2)\Phi^2+S\big)R=0.
}
Therefore,
\Eq{*}{
  S=-(4+3p)\Phi'-(3+5p+3p^2)\Phi^2.
}
Thus
\Eq{*}{
  \Psi_{f,g}=\frac{S+R}{2}
  &=\gamma\big|W_{f,g}^{1,0}\big|^p-\frac{4+3p}2\Phi'-\frac{3+5p+3p^2}2\Phi^2,\\
  \Psi_{F,G}=\frac{S-R}{2}
  &=-\gamma\big|W_{f,g}^{1,0}\big|^p-\frac{4+3p}2\Phi'-\frac{3+5p+3p^2}2\Phi^2,
}
which was to be proved.
\end{proof}

\Lem{N3}{
Let $\mu$ be a Borel probability measure on $[0,1]$ with $\mu_2>0=\mu_3$ and let $(f,g),(F,G)\in\C_6(I)$. If, for all $x\in I$,
\Eq{2+4+6d}{
m_{x;f,g;\mu}''(0)=m_{x;F,G;\mu}''(0), \qquad
m_{x;f,g;\mu}''''(0)=m_{x;F,G;\mu}''''(0) \qquad\mbox{and}\qquad
m_{x;f,g;\mu}''''''(0)=m_{x;F,G;\mu}''''''(0)
}
hold, then $\Phi_{f,g}=\Phi_{F,g}=:\Phi$ and either $\Psi_{f,g}=\Psi_{F,g}$ or there exists a nonzero real constant $\gamma$ such that \eq{R} holds with $p$ defined by \eq{p} and we have the following alternatives:
\begin{enumerate}[(i)]
 \item If $\mu_6=5\mu_2\mu_4$ and $\mu_4=3\mu_2^2$, then $\Phi$ is an at most first degree polynomial and $\Psi_{f,g}-\Psi_{F,g}$ is constant on $I$.
 \item If $\mu_6=5\mu_2\mu_4$ and $\mu_4\neq 3\mu_2^2$, then the equalites
 \Eq{*}{
    \Psi_{f,g}&=\gamma\big|W_{f,g}^{1,0}\big|^p
    -\frac{p+1}{3p}\frac{\Phi''}{\Phi}-\frac{2p+3}2\Phi'
    -\frac{2p^2+3p+4}6\Phi^2, \\
    \Psi_{F,G}&=-\gamma\big|W_{f,g}^{1,0}\big|^p
    -\frac{p+1}{3p}\frac{\Phi''}{\Phi}-\frac{2p+3}2\Phi'
    -\frac{2p^2+3p+4}6\Phi^2
 }
 hold on the subset $I\setminus\Phi^{-1}(0)=\{x\in I\mid \Phi(x)\neq0\}$.
 \item If $\mu_6\neq 5\mu_2\mu_4$ and $\mu_4=3\mu_2^2$, then there exists a real constant $\delta$ such that the equalities
\Eq{SS0}{
   \Psi_{f,g}=&\gamma
   +\delta\big|W_{f,g}^{1,0}\big|^{-1}-\frac{r}{2}\Phi'
   +\frac{r-5}{4} \Phi^2-\frac{3r-7}{12}\cdot\big|W_{f,g}^{1,0}\big|^{-1}
   \bigg(\int \Phi^3\big|W_{f,g}^{1,0}\big|\bigg), \\
   \Psi_{F,G}=&-\gamma
    +\delta\big|W_{f,g}^{1,0}\big|^{-1}-\frac{r}{2}\Phi'
   +\frac{r-5}{4} \Phi^2-\frac{3r-7}{12}\cdot\big|W_{f,g}^{1,0}\big|^{-1}
   \bigg(\int \Phi^3\big|W_{f,g}^{1,0}\big|\bigg)
}
hold on $I$, where 
\Eq{r}{
 r:=\dfrac{7\mu_6-45\mu_2^3}{3\mu_6-45\mu_2^3}.
}
 \item If $\mu_6\neq 5\mu_2\mu_4$ and $\mu_4\neq 3\mu_2^2$, then there exists a real constant $\delta$ such that the equalities
\Eq{SS}{
   \Psi_{f,g}=&\gamma\big|W_{f,g}^{1,0}\big|^p
   +\delta\big|W_{f,g}^{1,0}\big|^q+\frac{2(q-p)(p+1)-p+2}{6p}\Phi'
    +\frac{(q-p)(q+3p+1)(p+1)+p^2-2p}{6p} \Phi^2\\
   &+\frac{(q-p)(2p+q)(p+q+1)(p+1)}{6p}\cdot\big|W_{f,g}^{1,0}\big|^q
   \bigg(\int \Phi^3\big|W_{f,g}^{1,0}\big|^{-q}\bigg), \\
    \Psi_{F,G}=&-\gamma\big|W_{f,g}^{1,0}\big|^p
    +\delta\big|W_{f,g}^{1,0}\big|^q+\frac{2(q-p)(p+1)-p+2}{6p}\Phi'
    +\frac{(q-p)(q+3p+1)(p+1)+p^2-2p}{6p} \Phi^2\\
   &+\frac{(q-p)(2p+q)(p+q+1)(p+1)}{6p}\cdot\big|W_{f,g}^{1,0}\big|^q
   \bigg(\int \Phi^3\big|W_{f,g}^{1,0}\big|^{-q}\bigg)
}
hold on $I$, where 
\Eq{q}{
 q:=\frac{\mu_2}{\mu_4}\cdot\frac{10\mu_4^2-3\mu_6\mu_2-15\mu_4\mu_2^2}{\mu_6-5\mu_4\mu_2}.
}
\end{enumerate}
}

\begin{proof}
The condition $\mu_2>0$ implies that $\mu_4>0$ and $\mu_6>0$ are also valid. The equalities $\Phi_{f,g}=\Phi_{F,g}$ and \eq{R} are consequences of \lem{N1} and \lem{N2}, respectively. If $\gamma=0$, then $\Psi_{f,g}=\Psi_{F,g}$.
Therefore, in the rest of the proof, we may assume that $\gamma$ is not zero.

Using the formula for the sixth-order derivative by \prp{der}, the third condition of \eq{2+4+6d}, we arrive at
\Eq{*}{
&15\mu_2^3\big(8\Phi^3\Psi_{f,g}+6\Phi\Psi^2_{f,g}\big)
-15\mu_2\mu_4\big(3\Phi\Psi_{f,g}^2
      +\big(\Phi''+5\Phi'\Phi+4\Phi^3\big)\Psi_{f,g}+3\Phi^2\Psi'_{f,g}+\Phi\Psi''_{f,g}+2\Psi'_{f,g}\Psi_{f,g}\big) \\
    &+\mu_6\big(3\Phi\Psi^2_{f,g}+\big(7\Phi''+15\Phi'\Phi
      +4\Phi^3\big)\Psi_{f,g}+\big(12\Phi'+9\Phi^2\big)\Psi'_{f,g}+9\Phi\Psi''_{f,g}+6\Psi'_{f,g}\Psi_{f,g}+4\Psi'''_{f,g}\big)\\
    &=15\mu_2^3\big(8\Phi^3\Psi_{F,G}+6\Phi\Psi^2_{F,G}\big)
-15\mu_2\mu_4\big(3\Phi\Psi_{F,G}^2
      +\big(\Phi''+5\Phi'\Phi+4\Phi^3\big)\Psi_{F,G}+3\Phi^2\Psi'_{F,G}+\Phi\Psi''_{F,G}\\
      &+2\Psi'_{F,G}\Psi_{F,G}\big)
    +\mu_6\big(3\Phi\Psi^2_{F,G}+\big(7\Phi''+15\Phi'\Phi
      +4\Phi^3\big)\Psi_{F,G}+\big(12\Phi'+9\Phi^2\big)\Psi'_{F,G}+9\Phi\Psi''_{F,G}\\
      &+6\Psi'_{F,G}\Psi_{F,G}+4\Psi'''_{F,G}\big).
}
Thus, introducing $R:=\Psi_{f,g}-\Psi_{F,G}$ and $S:=\Psi_{f,g}+\Psi_{F,G}$, we can rewrite this equation as 
\Eq{RS}{
4\mu_6&R'''+(9\mu_6-15\mu_2\mu_4)\Phi R''+(9\mu_6-45\mu_2\mu_4)\Phi^2R'+12\mu_6\Phi'R'+(3\mu_6-15\mu_2\mu_4)(SR'+RS')\\&+(7\mu_6-15\mu_2\mu_4)\Phi''R+(15\mu_6-75\mu_2\mu_4)\Phi'\Phi R+(4\mu_6-60\mu_2\mu_4+120\mu_2^3)\Phi^3R\\&+(3\mu_6-45\mu_2\mu_4+90\mu_2^3)\Phi RS=0.
}
On the other hand, using \eq{R'}, \eq{R''}, and $R'''=(p\Phi''+3p^2\Phi'\Phi+p^3\Phi^3)R$, we obtain
\Eq{*}{
R\Big(&4\mu_6(p\Phi''+3p^2\Phi'\Phi+p^3\Phi^3)+(9\mu_6-15\mu_2\mu_4)(p\Phi'\Phi+p^2\Phi^3)+p(9\mu_6-45\mu_2\mu_4)\Phi^3\\
&+12p\mu_6\Phi'\Phi+(3\mu_6-15\mu_2\mu_4)(p\Phi S+S')
+(7\mu_6-15\mu_2\mu_4)\Phi''+(15\mu_6-75\mu_2\mu_4)\Phi'\Phi\\
&+(4\mu_6-60\mu_2\mu_4+120\mu_2^3)\Phi^3
+(3\mu_6-45\mu_2\mu_4+90\mu_2^3)\Phi S\Big)=0.
}
Due to $\gamma\neq0$, the function $R$ is nowhere zero and we get the following equation for $S$ and $\Phi$:
\Eq{EEE}{
3(&\mu_6-5\mu_2\mu_4)S'+(3\mu_6(p+1)-15\mu_2\mu_4(p+3)+90\mu_2^3)\Phi S+(\mu_6(4p+7)-15\mu_2\mu_4)\Phi''\\
&+(3\mu_6(4p^2+7p+5)-15\mu_2\mu_4(p+5))\Phi'\Phi+(\mu_6(4p^3+9p^2+9p+4)\\
&-15\mu_2\mu_4(p^2+3p+4)+120\mu_2^3)\Phi^3=0.
}

Consider first the case when $\mu_4=3\mu_2^2$. Then $p=0$ and \eq{EEE} simplifies to
\Eq{EE0}{
  3(\mu_6-15\mu_2^3)S'+3(\mu_6-15\mu_2^3)\Phi S+(7\mu_6-45\mu_2^3)\Phi''
  +15(\mu_6-15\mu_2^3)\Phi'\Phi+4(\mu_6-15\mu_2^3)\Phi^3=0.
}
If additionally $\mu_6=5\mu_2\mu_4=15\mu_2^3$, then 
\Eq{*}{
  60\mu_2^3\Phi''=0,
}
which yields that $\Phi$ is an at most first degree polynomial. This, together with the result of \lem{N2}, completes the proof of assertion (i). On the other hand, if $\mu_6\neq 5\mu_2\mu_4=15\mu_2^3$, then dividing the equation \eq{EE0} by $3(\mu_6-15\mu_2^3)\neq0$, we get
\Eq{*}{
  S'+\Phi S+r\Phi''+5\Phi'\Phi+\frac{4}{3}\Phi^3=0,
}
where $r$ is defined in \eq{r}.
This is a first-order inhomogeneous linear differential equation for $S$, whose general solution is of the form
\Eq{*}{
  S=2\delta\big|W_{f,g}^{1,0}\big|^{-1}-r\Phi'+\frac{r-5}{2} \Phi^2
   -\frac{3r-7}{6}\cdot\big|W_{f,g}^{1,0}\big|^{-1}
   \bigg(\int \Phi^3\big|W_{f,g}^{1,0}\big|\bigg),
}
where $\delta$ is an arbitrary real constant. Using the equalities
\Eq{Psi}{
  \Psi_{f,g}=\frac{S+R}{2} \qquad\mbox{and}\qquad \Psi_{F,G}=\frac{S-R}{2}
}
and \eq{R}, assertion (iii) follows directly.

Now consider the case when $\mu_4\neq 3\mu_2^2$ which is equivalent to $p\neq0$. 
If $\mu_6=5\mu_2\mu_4$, then the coefficient of $S'$ in the equality \eq{EEE} is zero. Using the definition \eq{p} of $p$ and the equality $\mu_6=5\mu_2\mu_4$, \eq{EEE} can be rewritten as
\Eq{*}{
  30\mu_2\mu_4p \Phi S+20\mu_2\mu_4(p+1)\Phi''
  +30\mu_2\mu_4(2p^2+3p)\Phi'\Phi+10\mu_2\mu_4(2p^3+3p^2+4p)\Phi^3=0.
}
This implies that
\Eq{*}{
  p \Phi S+20(p+1)\Phi''
  +30(2p^2+3p)\Phi'\Phi+10(2p^3+3p^2+4p)\Phi^3=0.
}
Therefore, using \eq{Psi} and \eq{R}, the second assertion results.

If $\mu_6\neq 5\mu_2\mu_4$, then a simple computation shows that
\Eq{*}{
  p+q+1=\frac{10\mu_2\mu_4}{5\mu_2\mu_4-\mu_6}p,
}
therefore we have that $p+q+1\neq0$ and the coefficient of $S'$ in the equality \eq{EEE} is also not zero. From the definitions of $p$ and $q$ it follows that
\Eq{*}{
\mu_4=\dfrac{3\mu_2^2}{p+1} \qquad\mbox{and}\qquad 
\mu_6=\dfrac{15\mu_2^3(q-p+1)}{(p+1)(p+q+1)}, 
}
respectively. Substituting these expressions into \eq{EEE} and then multiplying it by $-\dfrac{(p+1)(p+q+1)}{90\mu_2^3p}$, we arrive at
\Eq{NH}{
S'&-q\Phi S
+\frac{2p^2-2pq+3p-2q-2}{3p}\Phi''
+(2p^2-2pq+2p-3q+2)\Phi'\Phi\\
&+\frac13(2p^3-2p^2q+4p^2-3pq+2p-4q)\Phi^3=0.
}
This is an inhomogeneous first-order linear differential equation for $S$, whose general solution is of the following form
\Eq{*}{
  S=&2\delta\big|W_{f,g}^{1,0}\big|^q+\frac{2(q-p)(p+1)-p+2}{3p}\Phi'
    +\frac{(q-p)(q+3p+1)(p+1)+p^2-2p}{3p} \Phi^2\\
   &+\frac{(q-p)(2p+q)(p+q+1)(p+1)}{3p}\cdot\big|W_{f,g}^{1,0}\big|^q
   \bigg(\int \Phi^3\big|W_{f,g}^{1,0}\big|^{-q}\bigg).
}
This equality combined with \eq{R} and \eq{Psi} completes the proof of the last assertion of the lemma.
\end{proof}

\Cor{N3}{
Let $\mu$ be a Borel probability measure on $[0,1]$ with $\mu_2>0=\mu_3$ and $6\mu_6\mu_2^2-\mu_6\mu_4-5\mu_4^2\mu_2=0$ and let $(f,g),(F,G)\in\C_6(I)$. If, for all $x\in I$, the equalities in \eq{2+4+6d} are satisfied, then $\Phi_{f,g}=\Phi_{F,g}=:\Phi$ and we have the following alternatives
\begin{enumerate}[(i)]
 \item If $\mu_4=3\mu_2^2$, then $p=0$, $\Phi$ is an at most first degree polynomial and $\Psi_{f,g}-\Psi_{F,G}$ is a constant.
 \item If $\mu_4\neq3\mu_2^2$, then $p\neq0$ and there exist real constants $\alpha$ and $\beta$ such that the equalities
\Eq{SS+}{
   \Psi_{f,g}=\alpha\big|W_{f,g}^{1,0}\big|^p
   +\frac{2-p}{6p}\Phi'+\frac{p-2}{6}\Phi^2 \qquad\mbox{and}\qquad
    \Psi_{F,G}=\beta\big|W_{f,g}^{1,0}\big|^p
    +\frac{2-p}{6p}\Phi'+\frac{p-2}{6} \Phi^2
}
hold on $I$, where $p$ is given by \eq{p}.
\end{enumerate}
}

\begin{proof} The equality $\Phi_{f,g}=\Phi_{F,g}$ is a consequence of \lem{N1}.

If $\mu_4=3\mu_2^2$, then the moment condition $6\mu_6\mu_2^2-\mu_6\mu_4-5\mu_4^2\mu_2=0$ implies that $\mu_6=5\mu_4\mu_2$. Hence, we are in the alternative (i) of \lem{N3}, which yields that $\Phi$ is an at most first degree polynomial and $\Psi_{f,g}-\Psi_{F,G}$ is a constant on $I$. 

If $\mu_4\neq3\mu_2^2$, then $p\neq0$ and the moment condition $6\mu_6\mu_2^2-\mu_6\mu_4-5\mu_4^2\mu_2=0$ implies that $\mu_6\neq5\mu_4\mu_2$ and $p=q$, where $q$ is defined by \eq{q}. Now we satisfy the conditions of the alternative (iv) of \lem{N3}, hence \eq{SS+} is valid with $\alpha:=\gamma+\delta$ and $\beta:=\delta-\gamma$.
\end{proof}

As the main applications the above corollary, we restate and reprove the solution of the equality problems related to Bajraktarevi\'c and Cauchy means in the following two subsections.

For a real parameter $t\in\R$, introduce the sine and cosine type functions $S_t,C_t:\R\to\R$ by
\Eq{*}{
  S_t(x):=\begin{cases}
           \sin(\sqrt{-t}x) & \mbox{ if } t<0, \\
           x & \mbox{ if } t=0, \\
           \sinh(\sqrt{t}x) & \mbox{ if } t>0, 
         \end{cases}\qquad\mbox{and}\qquad
  C_t(x):=\begin{cases}
           \cos(\sqrt{-t}x) & \mbox{ if } t<0, \\
           1 & \mbox{ if } t=0, \\
           \cosh(\sqrt{t}x) & \mbox{ if } t>0. 
        \end{cases}
}
It is easily seen that the functions $S_t$ and $C_t$ form the fundamental system of solutions for the second-order homogeneous linear differential equation $h''=th$.

\subsection{\bf Equality of Bajraktarevi\'c means}

The first main result of this section is a rephrased form of the result of Losonczi \cite{Los99,Los06b} which characterized the equality of Bajraktarevi\'c means. In these papers Losonczi established 1+32 cases for the equality of these means. To deduce the result of Losonczi from the theorem below, the best is to elaborate condition (vi) where, beyond the canonical case (that is the equivalence of the generating functions), the equality is described in terms of two polynomials of at most second degree. In the subcases when, independently, these polynomials are constants, of first degree, of second degree with no, or with one or with two real roots, we can distinguish $6\times6=36$ subcases which then reduce to the cases considered by Losonczi.

\Thm{EBM}{Let $(f,g),(F,G)\in\C_6(I)$ and let $\mu:=\frac{\delta_0+\delta_1}{2}$. Then the following assertions are equivalent:
\begin{enumerate}[(i)]
 \item The means $M_{f,g;\mu}$ and $M_{F,G;\mu}$ are equal on $I^2$.
 \item The means $M_{f,g;\mu}$ and $M_{F,G;\mu}$ are equal near the diagonal of  $I^2$.
 \item For all $x\in I$, the equalities in \eq{2+4+6d} are satisfied.
 \item $\Phi_{f,g}=\Phi_{F,G}$ holds on $I$ and there exist real constants $\alpha$ and $\beta$ such that
 \Eq{*}{
    \Psi_{f,g}=\alpha\big(W_{f,g}^{1,0}\big)^2 \qquad\mbox{and}\qquad
    \Psi_{F,G}=\beta\big(W_{f,g}^{1,0}\big)^2.
 }
 hold on $I$.
 \item Either $(f,g)\sim(F,G)$ or there exist real constants $a,b,c,A,B,C,\gamma$ such that
 \Eq{Q}{
   af^2+bfg+cg^2=1 \qquad\mbox{and}\qquad AF^2+BFG+CG^2=1
 }
 and $W_{F,G}^{1,0}=\gamma W_{f,g}^{1,0}$.
 \item Either $(f,g)\sim(F,G)$ or there exist two real polynomials $P$ and $Q$ of at most second degree which are positive on the range of $f/g$ and $F/G$, respectively, and there exist real constants $\gamma$ and $\delta$ such that
 \Eq{PQ}{
   g=\frac{1}{\sqrt{P}}\circ\frac{f}{g}, \qquad
   G=\frac{1}{\sqrt{Q}}\circ\frac{F}{G}, \qquad\mbox{and}\qquad
   \bigg(\int\frac{1}{Q}\bigg)\circ\frac{F}{G}
   =\gamma\bigg(\int\frac{1}{P}\bigg)\circ\frac{f}{g}+\delta.
 }
 \item Either $(f,g)\sim(F,G)$ or there exist a strictly monotone function $\varphi:I\to\R$ and real constants $\alpha$ and $\beta$ such that
 \Eq{*}{
   (f,g)\sim(S_\alpha\circ\varphi,C_\alpha\circ\varphi) \qquad\mbox{and}\qquad
   (F,G)\sim(S_\beta\circ\varphi,C_\beta\circ\varphi)
 }
 \item Either $(f,g)\sim(F,G)$ or $M_{f,g;\mu}=A_\varphi=M_{F,G;\mu}$ holds on $I^2$ holds with $\varphi:=\int W_{f,g}^{1,0}$.
 \item Either $(f,g)\sim(F,G)$ or there exists a strictly monotone function $\varphi:I\to\R$ such that $M_{f,g;\mu}=A_\varphi=M_{F,G;\mu}$ holds on $I^2$.
\end{enumerate}}

\begin{proof}
The implication (i)$\Rightarrow$(ii) is obvious. The implication (ii)$\Rightarrow$(iii) is a direct consequence of \lem{N0}. Assume now that assertion (iii) is valid. Using that $\mu$ is of the form $\dfrac{\delta_0+\delta_1}{2}$, an easy computation shows that
\Eq{*}{
\muhat_1=\frac12 \qquad\mbox{and}\qquad
  \mu_n=\begin{cases}
        0 &\mbox{if } n\ \mbox{is\ odd}\\
         \frac{1}{2^n}  &\mbox{if }n\ \mbox{is\ even}
        \end{cases} \qquad(n\in\N\cup\{0\}).
}
Therefore, conditions $\mu_2>0=\mu_3$, $\mu_4\neq3\mu_2^2$ and $6\mu_6\mu_2^2-\mu_6\mu_4-5\mu_4^2\mu_2=0$ hold, whence using \eq{p}, we get $p=2$. Now the second alternative of \cor{N3} is applicable and it implies assertion (iv). 

To prove the implication (iv)$\Rightarrow$(v), assume that (iv) holds for some constants $\alpha,\beta\in\R$. If $\alpha=\beta$, then \lem{3} implies that $(f,g)\sim(F,G)$. Now consider the case when $\alpha\neq\beta$. The existence of some $\gamma$ such that the identity $W_{F,G}^{1,0}=\gamma W_{f,g}^{1,0}$ be valid is a direct consequence of the integration of the equality $\Phi_{f,g}=\Phi_{F,G}$. Applying implication (iv)$\Rightarrow$(ii) of \cite[Theorem 10]{PalZak20b}, we conclude that there exist real constants $a,b,c,A,B,C$ such that the equalities \eq{Q} hold.  Therefore, assertion (v) is valid. 

Assume now that (v) holds for some constants $a,b,c,A,B,C,\gamma$ and define
\Eq{*}{
  P(t):=at^2+bt+c \qquad\mbox{and}\qquad Q(t):=At^2+Bt+C\qquad(t\in\R).
}
Then, dividing the equalities in \eq{Q} side by side by $g^2$ and by $G^2$, we obtain that
\Eq{PQ+}{
  P\circ\frac{f}{g}=\frac{1}{g^2} \qquad\mbox{and}\qquad
  Q\circ\frac{F}{G}=\frac{1}{G^2}.
}
Therefore, $P$ and $Q$ are positive on the codomain of $f/g$ and $F/G$, respectively, and the first two equalities in \eq{PQ} hold. Furthermore, using \eq{PQ+}, we have
that
\Eq{*}{
  W_{f,g}^{1,0}=f'g-g'f=g^2\cdot\bigg(\frac{f}{g}\bigg)'
  =\bigg(\frac{1}{P}\circ\frac{f}{g}\bigg)\cdot\bigg(\frac{f}{g}\bigg)'
  =\bigg(\bigg(\int\frac{1}{P}\bigg)\circ\frac{f}{g}\bigg)'
}
and similarly,
\Eq{*}{
  W_{F,G}^{1,0}=\bigg(\bigg(\int\frac{1}{Q}\bigg)\circ\frac{F}{G}\bigg)'.
}
Applying the equality $W_{F,G}^{1,0}=\gamma W_{f,g}^{1,0}$, after integration we obtain that the third equality in \eq{PQ} is also valid for some real constant $\delta$. This shows that assertion (v) implies (vi).

Reversing the steps of the previous argument, one can easily see that assertion (vi) also implies (v), where $a,b,c$ and $A,B,C$ are the coefficients of the polynomials $P$ and $Q$, respectively.

To prove the implication (v)$\Rightarrow$(viii), assume that (v) holds for some constants $a,b,c,A,B,C,\gamma\in\R$. The implication (ii)$\Rightarrow$(iii) of \cite[Theorem 10]{PalZak20b} implies that 
\Eq{*}{
M_{f,g;\mu}=A_\varphi \qquad\mbox{and}\qquad M_{F,G;\mu}=A_\psi
}
hold on $I^2$ with $\varphi=\int W_{f,g}^{1,0}$ and $\psi=\int W_{F,G}^{1,0}$. Thus, using that $W_{f,g}^{1,0}=\gamma W_{F,G}^{1,0}$, we get $\varphi=\gamma\psi$. This result implies that $A_\varphi=A_\psi$ is satisfied on $I^2$. Therefore, assertion (viii) is valid.
 
The implications (viii)$\Rightarrow$(ix) and (ix)$\Rightarrow$(i) are obvious. Finally, the equivalence of the assertions (vii) and (ix) is a consequence of \cite[Corollary 9]{PalZak20b}.
\end{proof}

\subsection{\bf Equality of Cauchy means}
The second main result of this section is a rephrased form of the results of Losonczi \cite{Los00a,Los03a} which characterized the equality of Cauchy means and established 1+32 cases for the equality of these means. The results of Losonczi can easily be deduced from condition (vi) of the next theorem where, beyond the canonical case the equality is described in terms of two polynomials of at most second degree. Considering the same subcases as for \thm{EBM}, one can again distinguish $6\times6=36$ subcases which then reduce to the cases considered by Losonczi.

\Thm{ECM}{Let $(f,g),(F,G)\in\C_6(I)$ and let $\mu$ denote the Lebesgue measure restricted to $[0,1]$. Then the following assertions are equivalent:
\begin{enumerate}[(i)]
 \item The means $M_{f,g;\mu}$ and $M_{F,G;\mu}$ are equal on $I^2$.
 \item The means $M_{f,g;\mu}$ and $M_{F,G;\mu}$ are equal near the diagonal of  $I^2$.
 \item For all $x\in I$, the equalities in \eq{2+4+6d} are satisfied.
 \item $\Phi_{f,g}=\Phi_{F,G}$ and there exist constants $\alpha,\beta\in\R$ such that
 \Eq{PP}{
    \Psi_{f,g}=\alpha\big(W_{f,g}^{1,0}\big)^{\frac23} 
    +\frac13\Phi'-\frac29\Phi^2\qquad\mbox{and}\qquad
    \Psi_{F,G}=\beta\big(W_{f,g}^{1,0}\big)^{\frac23}+\frac13\Phi'-\frac29\Phi^2.
 }
 hold on $I$.
 \item Either $(f,g)\sim(F,G)$ or there exist real constants $a,b,c,A,B,C,\gamma$ such that
 \Eq{*}{
   af^2+bfg+cg^2=\big(W_{f,g}^{1,0}\big)^{\frac23} \qquad\mbox{and}\qquad 
   AF^2+BFG+CG^2=\big(W_{F,G}^{1,0}\big)^{\frac23}
 }
 and $W_{F,G}^{1,0}=\gamma W_{f,g}^{1,0}$.
 \item Either $(f,g)\sim(F,G)$ or there exist two real polynomials $P$ and $Q$ of at most second degree which are positive on the range of $f/g$ and $F/G$, respectively, and there exist real constants $\gamma$ and $\delta$ such that
 \Eq{PQC}{
   g=\bigg(\frac{f}{g}\bigg)'\bigg(\frac{1}{\sqrt{P^3}}\circ\frac{f}{g}\bigg), \quad
   G=\bigg(\frac{F}{G}\bigg)'\bigg(\frac{1}{\sqrt{Q^3}}\circ\frac{F}{G}\bigg), \quad\mbox{and}\quad
   \bigg(\int\frac{1}{Q}\bigg)\circ\frac{F}{G}
   =\gamma^{\frac13}\bigg(\int\frac{1}{P}\bigg)\circ\frac{f}{g}+\delta.
 }
 \item Either $(f,g)\sim(F,G)$ or there exist a strictly monotone differentiable function $\varphi:I\to\R$ and real constants $\alpha$ and $\beta$ such that
 \Eq{*}{
   (f,g)\sim(\varphi'\cdot S_\alpha\circ\varphi,\varphi'\cdot C_\alpha\circ\varphi) \qquad\mbox{and}\qquad
   (F,G)\sim(\varphi'\cdot S_\beta\circ\varphi,\varphi'\cdot C_\beta\circ\varphi)
 }
 \item Either $(f,g)\sim(F,G)$ or $M_{f,g;\mu}=A_\varphi=M_{F,G;\mu}$ holds on $I^2$ holds with $\varphi:=\int \big(W_{f,g}^{1,0}\big)^{\frac13}$.
 \item Either $(f,g)\sim(F,G)$ or there exists a strictly monotone function $\varphi:I\to\R$ such that $M_{f,g;\mu}=A_\varphi=M_{F,G;\mu}$ holds on $I^2$.
\end{enumerate}
}

\begin{proof}
The implication (i)$\Rightarrow$(ii) is clear. The implication (ii)$\Rightarrow$(iii) is a direct consequence of \lem{N0}. Assume now that assertion (iii) is valid. Using that $\mu$ is the Lebesgue measure restricted to $[0,1]$, it is easily seen that
\Eq{*}{
\muhat_1=\frac12 \qquad\mbox{and}\qquad
  \mu_n=\begin{cases}
        0 &\mbox{if } n\ \mbox{is\ odd}\\
         \dfrac{1}{(n+1)2^n}  &\mbox{if }n\ \mbox{is\ even}
        \end{cases} \qquad(n\in\N\cup\{0\}).
}
Consequently, conditions $\mu_2>0=\mu_3$, $\mu_4\neq3\mu_2^2$ and $6\mu_6\mu_2^2-\mu_6\mu_4-5\mu_4^2\mu_2=0$ are valid and
using \eq{p}, we get $p=\frac23$. Hence the second alternative of \cor{N3} is applicable. Thus $\Phi_{f,g}=\Phi_{F,G}$ holds and the equalities in \eq{SS+} reduce to \eq{PP}, which completes the proof of assertion (iv). 

To prove the implication (iv)$\Rightarrow$(v), assume that (iv) holds for some constants $\alpha,\beta\in\R$. If $\alpha=\beta$, then \lem{3} implies that $(f,g)\sim(F,G)$. Now consider the case when $\alpha\neq\beta$. The existence of some $\gamma$ such that the identity $W_{F,G}^{1,0}=\gamma W_{f,g}^{1,0}$ be valid is a direct consequence of the integration of the equality $\Phi_{f,g}=\Phi_{F,G}$. The equalities in \eq{PP} show that the expressions
\Eq{exp}{
 \frac{3W^{3,0}_{f,g}+12W^{2,1}_{f,g}}{\Big(W^{1,0}_{f,g}\Big)^{\frac53}}
 -5\frac{\Big(W^{2,0}_{f,g}\Big)^2}{\Big(W^{1,0}_{f,g}\Big)^{\frac83}}\qquad\mbox{and}\qquad
 \frac{3W^{3,0}_{F,G}+12W^{2,1}_{F,G}}{\Big(W^{1,0}_{F,G}\Big)^{\frac53}}
 -5\frac{\Big(W^{2,0}_{F,G}\Big)^2}{\Big(W^{1,0}_{F,G}\Big)^{\frac83}}
}
are constants. Therefore, applying implication (vii)$\Rightarrow$(iv) of \cite[Theorem 7]{LovPalZak20}, we conclude that there exist real constants $a,b,c,A,B,C$ that validate assertion our (v). 

Assume now that (v) holds for some constants $a,b,c,A,B,C,\gamma$ and define
\Eq{*}{
  P(t):=at^2+bt+c \qquad\mbox{and}\qquad Q(t):=At^2+Bt+C\qquad(t\in\R).
}
Then, dividing the equalities in \eq{Q} side by side by $g^2$ and by $G^2$, we obtain that
\Eq{PQ++}{
  P\circ\frac{f}{g}=\frac{\Big(W^{1,0}_{f,g}\Big)^{\frac23}}{g^2} \qquad\mbox{and}\qquad
  Q\circ\frac{F}{G}=\frac{\Big(W^{1,0}_{F,G}\Big)^{\frac23}}{G^2}.
}
Therefore, $P$ and $Q$ are positive on the codomain of $f/g$ and $F/G$, respectively. Using the identities $W^{1,0}_{f,g}=g^2(f/g)'$ and $W^{1,0}_{F,G}=G^2(F/G)'$, the above equalities yield the first two equations in \eq{PQC}. Furthermore, using \eq{PQ++}, we have
that
\Eq{*}{
  \bigg(\bigg(\int\frac{1}{P}\bigg)\circ\frac{f}{g}\bigg)'
  =\bigg(\frac{1}{P}\circ\frac{f}{g}\bigg)\cdot\bigg(\frac{f}{g}\bigg)'
  =\frac{g^2}{\Big(W^{1,0}_{f,g}\Big)^{\frac23}}\cdot\frac{W^{1,0}_{f,g}}{g^2}
  =\Big(W_{f,g}^{1,0}\Big)^{\frac13}
}
and similarly,
\Eq{*}{
  \bigg(\bigg(\int\frac{1}{Q}\bigg)\circ\frac{F}{G}\bigg)'
  =\Big(W_{F,G}^{1,0}\Big)^{\frac13}.
}
Applying the equality $W_{F,G}^{1,0}=\gamma W_{f,g}^{1,0}$, after integration we obtain that the third equality in \eq{PQC} is also valid for some real constant $\delta$. This completes the proof of the implication (v)$\Rightarrow$(vi).

Reversing the steps of the previous argument, one can easily see that assertion (vi) also implies (v), where $a,b,c$ and $A,B,C$ are the coefficients of the polynomials $P$ and $Q$, respectively.

To prove the implication (v)$\Rightarrow$(viii), assume that (v) holds for some constants $a,b,c,A,B,C,\gamma\in\R$. The implication (iv)$\Rightarrow$(vi) of \cite[Theorem 7]{LovPalZak20} implies that 
\Eq{*}{
M_{f,g;\mu}=A_\varphi \qquad\mbox{and}\qquad M_{F,G;\mu}=A_\psi
}
hold on $I^2$ with $\varphi=\int \Big(W_{f,g}^{1,0}\Big)^{\frac13}$ and $\psi=\int \Big(W_{F,G}^{1,0}\Big)^{\frac13}$. Thus, using that $W_{f,g}^{1,0}=\gamma W_{F,G}^{1,0}$, we get $\varphi=\gamma^{\frac13}\psi$. This equality yields that $A_\varphi=A_\psi$ is satisfied on $I^2$. Therefore, assertion (viii) is valid.
 
The implications (viii)$\Rightarrow$(ix) and (ix)$\Rightarrow$(i) are obvious. Finally, the equivalence of the assertions (vii) and (ix) is a consequence of the main result of the paper \cite{KisPal19}.
\end{proof}

\subsection{Conclusion and open problems}

We have to stress that the different assertions of \thm{EBM} and \thm{ECM} require different order of regularity. Obviously, assertions (i), (ii), (vi), (vii), (ix) make sense in the regularity class $\C_0(I)$. For (v) and (viii) one has to take the unknown functions from $\C_1(I)$. Finally, assertions (iv) and (iii) require the regularity class $\C_2(I)$ and $\C_6(I)$, respectively. \\\indent
One can also see that some of the implications described in the above proof are valid with smaller order regularity assumptions. For instance, (i) implies (ii) and (ix) implies (i) in the class $\C_0(I)$. For the implications (ii)$\Rightarrow$(iii)$\Rightarrow$(iv), we need $\C_6(I)$. The proof of the implication (iv)$\Rightarrow$(v) requires $\C_2(I)$, while the equivalence of assertions (v) and (vi), and the implications (v)$\Rightarrow$(viii)$\Rightarrow$(ix) can be verified in the regularity class $\C_1(I)$.

Based on the above observations, we can formulate the following three open problems.
\begin{enumerate}
 \item Provided that $(f,g)$, $(F,G)\in\C_0(I)$ and $\mu:=\frac{\delta_0+\delta_1}{2}$ (resp.\ $\mu$ is the Lebesgue measure on $[0,1]$), are assertions (i), (ii), (vi), (vii), and (ix) of \thm{EBM} (resp.\ \thm{ECM}) equivalent to each other?
\item Provided that $(f,g)$, $(F,G)\in\C_1(I)$ and $\mu:=\frac{\delta_0+\delta_1}{2}$ (resp.\ $\mu$ is the Lebesgue measure on $[0,1]$), are assertions (i), (ii), (v), (vi), (vii), (viii), and (ix) of \thm{EBM} (resp.\ \thm{ECM}) equivalent to each other?
\item Provided that $(f,g)$, $(F,G)\in\C_2(I)$ and $\mu:=\frac{\delta_0+\delta_1}{2}$ (resp.\ $\mu$ is the Lebesgue measure on $[0,1]$), are assertions (i), (ii), (iv), (v), (vi), (vii), (viii), and (ix) of \thm{EBM} (resp.\ \thm{ECM}) equivalent to each other?
\end{enumerate}


\end{document}